\documentclass[12pt,a4paper]{article}
\usepackage{amsfonts}

\usepackage{amssymb}

\usepackage{amsmath}

\newtheorem{theorem}{Theorem}[section]

\newtheorem{lemma}[theorem]{Lemma}

\date{}

\begin{document}

\title{ Gr\"{o}bner-Shirshov Bases for Associative Algebras with
Multiple Operators and Free Rota-Baxter Algebras\footnote{Dedicated
to Professor Yu.G. Reshetnyak on the occasion of his 80th birthday.}
\footnote{Supported by the NNSF of China (No.10771077) and the NSF
of Guangdong Province (No.06025062).} }

\author{
L. A. Bokut\footnote {Supported by RFBR 01-09-00157, LSS--344.2008.1
and SB RAS Integration grant No. 2009.97 (Russia).} \\
{\small \ School of Mathematical Sciences, South China Normal
University}\\
{\small Guangzhou 510631, P. R. China}\\
{\small Sobolev Institute of Mathematics, Russian Academy of
Sciences}\\
{\small Siberian Branch, Novosibirsk 630090, Russia}\\
{\small bokut@math.nsc.ru}\\
\\
 Yuqun
Chen\\
{\small \ School of Mathematical Sciences, South China Normal
University}\\
{\small Guangzhou 510631, P. R. China}\\
{\small yqchen@scnu.edu.cn}\\
\\
Jianjun Qiu\footnote {Corresponding author.}\\
{\small \ Mathematics and Computational  Science School }\\
{\small \ Zhanjiang  Normal University}\\
{\small \ Zhanjiang  524048, China}\\
{\small \ jianjunqiu@126.com}} \vspace{4mm}

\maketitle \noindent\textbf{Abstract:} In this paper, we establish
the Composition-Diamond lemma for associative algebras with multiple
linear operators.  As applications, we obtain Gr\"{o}bner-Shirshov
bases of free Rota-Baxter algebra, $\lambda$-differential algebra
and $\lambda$-differential Rota-Baxter algebra, respectively.  In
particular,  linear bases of these three free algebras are
respectively obtained, which are essentially the same or similar the
recent results obtained to those obtained by K. Ebrahimi-Fard -- L.
Guo, and L. Guo -- W. Keigher by  using other methods.



\section{Introduction}
Gr\"{o}bner and Gr\"{o}bner-Shirshov bases theory were initiated
independently by A.I. Shirshov \cite{Sh62a, Sh62b}, H. Hironaka
\cite{Hi64} and B. Buchberger \cite{Bu65, Bu70}.

In \cite{Sh62a}, Shirshov first defined Gr\"{o}bner-Shirshov basis
of a (anti-) commutative non-associative algebra presentation with a
fix monomial well ordering under the name ``reduced set of
relations" and he proved the ``Composition-Diamond Lemma" (Lemma 2
of \cite{Sh62a}) for any reduced set of a free (anti-) commutative
algebra. Also Shirshov mentioned that his theory is equally valid
for a non-associative algebra presentation and a  free
non-associative algebra. Twelve years earlier, A.I. Zhukov
\cite{Zh50} invented a weaker forms of Gr\"{o}bner-Shirshov bases of
non-associative algebras and Composition-Diamond lemma for a free
non-associative algebra. The difference is that Zhukov \cite{Zh50},
contrary to Shirshov \cite{Sh62a}, did not use any linear order of
monomials. Instead, he defined the ``leading" word for a
non-associative polynomial as any of its monomial of maximal degree.

In \cite{Sh62b}, Shirshov invented a much more involved
Gr\"{o}bner-Shirshov bases theory for Lie algebras explicitly and
for associative algebras implicitly, and  he proved the
Composition-Diamond Lemma for a free Lie algebra again explicitly
and for a free associative algebra implicitly (Lemma 3 of
\cite{Sh62b}). The main differences with \cite{Sh62a} are as
follows. For Lie algebra case, one needs to use a very special
linear basis of a free Lie algebra, Lyndon--Shirshov basis (see
Shirshov (\cite{Sh58}, 1958), Chen--Fox--Lyndon (\cite{CFL58},
1958)). Since that time, the Lyndon--Shirshov basis and
Lyndon--Shirshov words (due to R. Lyndon (\cite{Lyn54}, 1954), and
independently to Shirshov (\cite{Sh58}, 1958))  play a significant
role in algebra and combinatorics (see, for example, \cite{Bah87,
BMPZ92, BoKu94, Coh65, Lo83, Re93}; they are often called Lyndon
basis and Lyndon words, see, for example, \cite{Lo83, Re93}). Also
in Lie and associative algebra cases, one needs to use a notion of
composition (of intersection) of two Lie (associative) polynomials
($s$-polynomial under the later Buchberger's terminology
\cite{Bu70}). Shirshov did it in two steps. First, he defined the
composition of two Lie polynomials as associative (non-commutative)
polynomials, and then he put in special bracketing using his the
Special Bracketing Lemma 4 of \cite{Sh58}. Finally for Lie algebra
case, one needs to use Shirshov's reduction algorithm to eliminate
the leading (in degree-lexicographical ordering) monomial of one Lie
polynomial in another Lie polynomial (again using  Lemma 4 of
\cite{Sh58}); for associative case, the reduction algorithm is just
non-commutative Euclidean algorithm. An ingenious Shirshov's
algorithm was defined for any reduced set $S$ of Lie (associative)
polynomials: to joint to $S$ a composition of two polynomials from
$S$ and to apply the Shirshov's (Euclidean) reduction algorithm to
the result set, and then by repeating this process again. A detailed
description of this algorithm for associative polynomials (more
pricesely, for associative binomials) can be found, for example, in
the text \cite{Eps92} under the name Knuth-Bendix algorithm (cf. D.
Knuth, P. Bendix (\cite{KnuBen70}, 1970)). Shirshov algorithm is in
general infinite as well as the Knuth--Bendix algorithm (in contrary
with Shirshov algorithm, the latter is applied to any universal
algebra presentation (not only to Lie, associative, or
associative-commutative algebra presentations, or a semigroup
presentation), it does not use any specific normal forms (in free
Lie, associative, associative--commutative algebras, or free
semigroups, free groups), and it does not use any effectively
defined ``composition" ($s$-polynomial); it means that Knuth-Bendix
algorithm is too slowly to be effective for groups, semigroups, Lie
algebras, associative algebras, associative-commutative algebras,
and so on). Using transfinite induction, for any reduced set $S$ of
Lie (associative) polynomials, one may find a reduced set $S^c$ such
that any composition of polynomials from $S^c$ is reduced to zero by
using Shirshov (Euclidean) reduction algorithm, and $S$ and $S^c$
generate the same ideal of the free Lie (associative) algebra. The
set $S^c$ is now called a Gr\"{o}bner-Shirshov basis for $S$ or a
Gr\"{o}bner-Shirshov basis of the ideal $Id(S)$ generated by $S$.
Shirhsov proved the Composition-Diamond Lemma for any subset $S$ of
a free Lie algebra (Lemma 3 of \cite{Sh62b}).

It is clear that Shirshov algorithm as well as the definition of
$S^c$ and the proof of Composition-Diamond Lemma for $S^c$ are
equally valid for associative polynomials. In this case, one dose
not need to use Lyndon--Shirshov words, Special Bracketing Lemma,
and Shirshov reduction  algorithm (just Euclidean algorithm for
non-commutative polynomials will be enough). An easy proof of the
Composition-Diamond Lemma for associative polynomials is just
following the proof of Lemma 3 of \cite{Sh62b} (cf. \cite{BoKu94}).
Later the Composition-Diamond Lemma for associative polynomials was
explicitly formulated and proved in a paper by L.A. Bokut,
\cite{Bo76}) as a specialization of the Shirshov's Lemma and also by
G. Bergman in \cite{Be78}).

In recent years there were quite a few results on
Gr\"{o}bner-Shirshov bases for associative algebras, semisimple Lie
(super) algebras, irreducible modules, Kac--Moode algebras, Coxeter
groups, braid groups, conformal algebras, free inverse semigroups,
Lodays' dialgebras, Leibniz algebras and so on, see, for example,
monographs \cite{BoKu94, MikZ95}, and papers \cite{Bo08, Bo09,
BoCS07, BoCL08, BoCZ08,  BoFK04, BoFKS03, BoKLM99, BoS01,  KL00a,
KL00b, KLLO02, KLLP07, La00a, La00b, La00c, La05,  Mik96, Ufn95,
VMik96}, and also surveys \cite{BoC08, BoC05, BoFKK00, BoK00,
BoK05}. Actually, conformal algebras and Loday's dialgebras are
multioperator algebras (conformal algebras have infinitely many
bilinear operations $x(n)y, n=0,1,2,...$ and a linear operation
$D(x)$, dialgebras have two bilinear operations, left and right
products).

 The concept of multioperator algebra ($\Omega$-algebra) was
invented by A.G. Kurosh \cite{Ku60} under an influence of the
concept of multioperator group of P.J . Higgins \cite{Hi56}. In
\cite{Ku60}, Kurosh generalized on a free $\Omega$-algebra  his
result \cite{Ku47} that any subalgebra of a free non-associative
 algebra is again free (the Nielsen--Schreier property). According to Kurosh, free
$\Omega$-algebras would enjoy many combinatorial properties of free
non-associative  algebras. Indeed, for example, recently V. Drensky
and R. Holtkamp \cite{DH08} have  proved an analogy of above
Shirshov--Zhukov's Composition-Diamond lemma for free
$\Omega$-algebras. It means that they extended Shirshov--Zhukov's
Gr\"{o}bner-Shirshov bases theory from non-associative algebras to
$\Omega$-algebras.  In 1969, Kurosh initiated to publish in the
Uspehi Mat. Nauk (Russian Math. Surveys) a series of papers by V.A.
Artamonov, M.S.  Burgin, F.I. Kizner, S.V. Polin, and Yu.K. Rebane
on the theory of $\Omega$-algebras and wrote a special introduction
paper on this series \cite{Ku69}. He predicted that the theory of
$\Omega$-algebras would have a role in future algebra as a subject
between the theory of (non-associative) rings and the theory of
universal algebras. Now after 40 years, his prediction is fulfilled.

  Multioperator Rota--Baxter algebras were invented by G. Baxter \cite{Bax60} and
  was given the current name after an important paper by G.-C. Rota \cite{Ro69}.

  Let $K$ be a commutative (associative) ring with unity and $\lambda\in K$.
  A Rota--Baxter algebra of weight $\lambda$ is an associative $K$-algebra $R$
  with a linear operator $P(x)$ satisfying the Rota--Baxter identity
$$
P(x)P(y)=P(xP(y)) + P(P(x)y) +\lambda P(xy), x,y\in R.
$$

A differential algebra of weight $\lambda$, called
$\lambda$-differential algebra, is an associative $K$-algebra with
linear $\lambda$-differential operator $D(x)$ such that
$$D(xy)=
D(x)y+xD(y)+\lambda D(x)D(y)$$
for any $x,y\in R$. L. Guo and W.
Keigher \cite{GK08} introduced a notion of $\lambda$-differential
Rota-Baxter algebra which is a Rota-Baxter $K$-algebra of weight
$\lambda$ with $\lambda$-differential operator $D$ such that
$DP=Id_R$.

There have been some constructions of free (commutative) Rota-Baxter
algebras. In this aspect, G.-C. Rota \cite{Ro69} and P. Cartier
\cite{Ca72} gave the explicit constructions of free commutative
Rota-Baxter algebras of weight $\lambda=1$, which they  called
shuffle Baxter and standard Baxter algebras, respectively. Recently
L. Guo and W. Keigher \cite{GK00a, GK00b} constructed the free
commutative Rota-Baxter algebras (with unit or without  unit) for
any $\lambda \in K$ using the mixable shuffle product. These are
called the mixable shuffle product algebras, which generalize the
classical construction of shuffle product algebras (see, for
example, \cite{Lo83}). K. Ebrahimi-Fard and L. Guo \cite{EG08a}
further constructed the free associative Rota-Baxter algebras by
using the Rota-Baxter words.

E.  Kolchin \cite{ko} considered  the differential algebras and
 constructed
a free   differential algebra. L. Guo and W. Keigher
 \cite{GK08} dealt with a $\lambda$-generalization  of this algebras.
 In \cite{GK08},  the free $\lambda$-differential Rota-Baxter
 algebra was  obtained
  by using the free Rota-Baxter algebra  on  planar
 decorated rooted trees.

K. Ebrahimi-Fard and  L. Guo  \cite{EG08b} use rooted trees and
forests to give explicit construction of free noncommutative
Rota-Baxter algebras on modules and sets. K. Ebrahimi-Fard, J. M.
Gracia-Bondia, and F. Patras \cite{EGP07} gave solution of the word
problem for free non-commutative Rota-Baxter algebra. A free
 Rota-Baxter  algebra is constructed on decorated trees by M. Aguiar
and M. Moreira \cite{AM}.

In this paper, we deal with  associative algebras with multiple
linear operators. We construct free associative algebras with
multiple linear operators and establish the Composition-Diamond
lemma for such algebras. As applications, we give the
Gr\"{o}bner-Shirshov bases of free Rota-Baxter algebra,
$\lambda$-differential algebra and $\lambda$-differential
Rota-Baxter algebra, respectively. Thus the linear bases of these
three free algebras can be hence obtained, respectively, which are
essentially the same or similar to those obtained  by using other
methods in the papers  by K. Ebrahimi-Fard and L. Guo \cite{EG08a}
and L. Guo  and W. Keigher \cite{GK08}.

Let $X$ be a set. Then we denote by $X^*$ and $S(X)$ the free monoid
and free semigroup generated by $X$, respectively. Let $X$ be a well
ordered set. By the deg-lex ordering on $X^*$, we mean the ordering
  which compares two words first by their degrees and then
compares lexicographically.

\section{ Free  associative algebra with multiple operators}

A free associative algebra with unary operators are constructed by
L. Guo \cite{Gu09}. A free associative algebra with one operator is
constructed by K. Ebrahimi-Fard and L. Guo \cite{EG08a}. In this
section, we construct the free  associative algebra with multiple
linear operators.

Let $K$ be a commutative ring with unit. An associative algebra with
multiple linear operators is an associative $K$-algebra $R$ with a
set
 $\Omega$ of multilinear operators (operations).

Let $X$ be a set and
$$
\Omega=\bigcup_{n=1}^{\infty}\Omega_{n},
 $$
where $\Omega_{n}$ is the set of $n$-ary operations, for example,
ary $(\delta)=n$ if $\delta\in \Omega_n$.

Define
$$\mathfrak{S}_{0}=S(X_0),\  X_0=X,$$
$$
\mathfrak{S}_{1}=S(X_1), \  X_1=X\cup \Omega(\mathfrak{S}_{0}) \\
$$
where
$$
\Omega(\mathfrak{S}_{0})=\bigcup\limits_{t=1}^{\infty}\{\delta(u_1,u_2,\dots,u_t)|\delta\in
\Omega_t, u_i\in \mathfrak{S}_{0} , \ i=1,2,\dots,t\}.
$$
For $n>1$,  define
$$
\mathfrak{S}_{n}=S(X_n), \ X_n=X\cup \Omega(\mathfrak{S}_{n-1})
$$
where
$$
\Omega(\mathfrak{S}_{n-1})=\bigcup\limits_{t=1}^{\infty}\{\delta(u_1,u_2,\dots,u_t)|\delta\in
\Omega_t, u_i\in \mathfrak{S}_{n-1} , \ i=1,2,\dots,t\}.
$$

Then we have
$$
\mathfrak{S}_{0}\subset\mathfrak{S}_{1}\subset\cdots
\subset\mathfrak{S}_{n}\subset\cdots.
$$
Let
$$
\mathfrak{S}(X)=\bigcup_{n\geq0}\mathfrak{S}_{n}.
$$
Then, it is easy to see that $\mathfrak{S}(X)$ is a semigroup such
that $ \Omega(\mathfrak{S}(X))\subseteq \mathfrak{S}(X). $

Now, we describe an element in $ \mathfrak{S}(X)$  by a labeled
reduced planar rooted forest, see \cite{DH08} and \cite{EG08b}.

For $x\in X$, we now describe  $x$ by a labeled reduced planar
rooted tree $\bullet_{x}$. Thus, for any  element $u=x_1 x_2\dots
x_n\in S(X)$, $u$ can be described by
$$
\bullet_{x_1}\sqcup \bullet_{x_2}\sqcup \dots \sqcup\bullet_{x_n}
$$
which is called a  labeled reduced planar rooted forest.  For any
$u\in \Omega(\mathfrak{S}_{0})$, say $u=\theta_n(u_1, \dots ,u_n)$
where each $u_i\in \mathfrak{S}_{0}$, $u$ can be described by

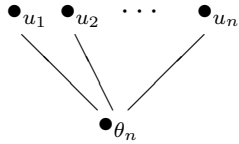
\setlength{\unitlength}{1cm}
\begin{picture}(4,2.0)
\put(7,0.5){\line(1,1){1}} \put(6.8,0.5){\line(-1,2){0.5}}
\put(6.6,0.5){\line(-1,1){1}}
\put(6.6,0.2){$\bullet_{\theta_n}$}\put(6.1,1.7){$\bullet_{u_2}$}
\put(5.4,1.7){$\bullet_{u_1}$}\put(7.9,1.7){$\bullet_{u_n}$}\
\put(6.9,1.7){$\cdots$}
\end{picture}\\
For example,  if $u=\theta_3(u_1, u_2,u_3)\in
\Omega(\mathfrak{S}_{0})$ where $u_1=x_1 x_2, \ u_2=x_3, \ u_3=x_4
x_5$, then

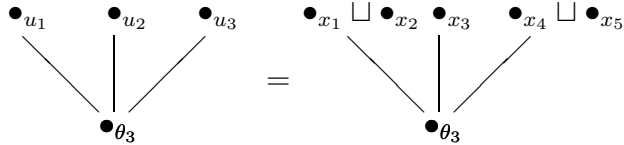
\setlength{\unitlength}{1cm}
\begin{picture}(4,2.0)
 \put(5,0.5){\line(1,1){1}} \put(4.8,0.5){\line(0,1){1}}
\put(4.6,0.5){\line(-1,1){1}} \put(4.6,0.2){$\bullet_{\theta_3}$}
\put(4.6,0.2){$\bullet_{\theta_3}$}\put(4.7,1.7){$\bullet_{u_2}$}
\put(3.4,1.7){$\bullet_{u_1}$}\put(5.9,1.7){$\bullet_{u_3}$}\put(6.8,
0.8){=}
\end{picture}
\setlength{\unitlength}{1cm}
\begin{picture}(4,2.0)
\put(5,0.5){\line(1,1){1}} \put(4.8,0.5){\line(0,1){1}}
\put(4.6,0.5){\line(-1,1){1}} \put(4.6,0.2){$\bullet_{\theta_3}$}
\put(4.6,0.2){$\bullet_{\theta_3}$}\put(4.7,1.7){$\bullet_{x_3}$}
\put(3,1.7){$\bullet_{x_1}\sqcup\bullet_{x_2}$}\put(5.7,1.7){$\bullet_{x_4}\sqcup\bullet_{x_5}$}
\end{picture}\\
Now, for any  $v=v_1v_2\cdots v_n\in \mathfrak{S}_{1}=S(X_1) \
(v_i\in X_1, \ i=1,\dots,n)$, $u$ can be described by labeled
reduced planar rooted forest:
$$
\bullet_{v_1}\sqcup\bullet_{v_2}\sqcup\cdots\sqcup\bullet_{v_n}
$$
For example, if $u=\theta_3(u_1, u_2,u_3)\theta_2(u_4,u_5)$ where
$u_1=x_1 x_2$, $u_2=x_3$, $u_3=x_4 x_5,\ u_4= x_6,\ u_5=x_7x_8$,
then $u$ can be described by \setlength{\unitlength}{1cm}

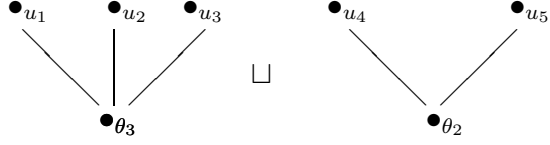
\begin{picture}(4,2.0)
\put(5,0.5){\line(1,1){1}} \put(4.8,0.5){\line(0,1){1}}
\put(4.6,0.5){\line(-1,1){1}} \put(4.6,0.2){$\bullet_{\theta_3}$}
\put(4.6,0.2){$\bullet_{\theta_3}$}\put(4.7,1.7){$\bullet_{u_2}$}
\put(3.4,1.7){$\bullet_{u_1}$}\put(5.7,1.7){$\bullet_{u_3}$}
\put(6.6, 0.8){$\sqcup$} \put(9.1 ,0.5){\line(1,1){1}}\put(8.9
,0.5){\line(-1,1){1}}\put(8.9,0.2){$\bullet_{\theta_2}$}\put(7.6
,1.7){$\bullet_{u_4}$}\put(10,1.7){$\bullet_{u_5}$}
\end{picture}\\
\noindent that is by

\setlength{\unitlength}{1cm}
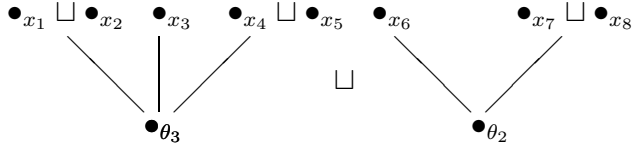
\begin{picture}(4,2.0)
\put(5,0.5){\line(1,1){1}} \put(4.8,0.5){\line(0,1){1}}
\put(4.6,0.5){\line(-1,1){1}} \put(4.6,0.2){$\bullet_{\theta_3}$}
\put(4.6,0.2){$\bullet_{\theta_3}$}\put(4.7,1.7){$\bullet_{x_3}$}
\put(2.8,1.7){$\bullet_{x_1}\sqcup\bullet_{x_2}$}\put(5.7,1.7){$\bullet_{x_4}\sqcup\bullet_{x_5}$}
\put(7.1, 0.8){$\sqcup$} \put(9.1 ,0.5){\line(1,1){1}}\put(8.9
,0.5){\line(-1,1){1}}\put(8.9,0.2){$\bullet_{\theta_2}$}\put(7.6
,1.7){$\bullet_{x_6}$}\put(9.5
,1.7){$\bullet_{x_7}\sqcup\bullet_{x_{8}}$}
\end{picture}\\
Continue this step, we can describe any element in $\mathfrak{S}(X)$
by labeled reduced planar rooted forest.

Therefore, with the above way, each element in $\mathfrak{S}(X)$
corresponds uniquely a labeled reduced planar rooted forest.

\ \

 For any $u\in \mathfrak{S}(X)$, $dep(u)=\mbox{min}\{n|u\in\mathfrak{S}_{n} \}$
 is called the depth of
 $u$.

Let   $K\langle X; \Omega\rangle$ be the $K$-algebra spanned by
$\mathfrak{S}(X)$.  Then, the element in $\mathfrak{S}(X)$ (resp.
$K\langle X; \Omega\rangle$) is called a $\Omega$-word (resp.
$\Omega$-polynomial). If $u\in X\cup \Omega(\mathfrak{S}(X))$, we
call $u$ a prime $\Omega$-word and define $bre(u)=1$ (the breadth of
$u$). If $u=u_1u_2\cdots u_n\in\mathfrak{S}(X)$, where $u_i$ is
prime $\Omega$-word for all $i$, then we define $bre(u)=n$.

Extend linearly each $\delta\in\Omega_n$,
$$\delta:\mathfrak{S}(X)^n\rightarrow \mathfrak{S}(X), \
(x_1,x_2,\cdots,x_n)\mapsto \delta(x_1,x_2,\cdots,x_n)
$$
to $K\langle X; \Omega\rangle$. Then, it is easy to see that
$K\langle X; \Omega\rangle$ is a free associative algebra with
multiple linear operators $\Omega$ on set $X$.

\section{ Composition-Diamond lemma for associative algebras with
multiple operators}

In this section, we introduce the notions of Gr\"{o}bner-Shirshov
bases  for associative algebras with multiple linear operators and
establish the Composition-Diamond lemma for such algebras.

Let $K\langle X;\Omega\rangle$ be free associative algebra with
multiple linear operators $\Omega$ on $X$ and $\star\notin X$. By a
$\star$-$\Omega$-word we mean any expression in $\mathfrak{S}(X\cup
\{\star\})$ with only one occurrence of $\star$. The set of all
$\star$-$\Omega$-words on $X$we define by $\mathfrak{S}^\star (X)$.

Let $u$ be a $\star$-$\Omega$-word and $s\in K\langle
X;\Omega\rangle$. Then we call
$$
u|_{s}=u|_{\star\mapsto s}
$$
an $s$-$\Omega$-word.

In other words, an $s$-$\Omega$-word $u|_{\star\mapsto s}$  means
that we have replaced the leaf $\star$ of $u$ by $s$.

For example, \\

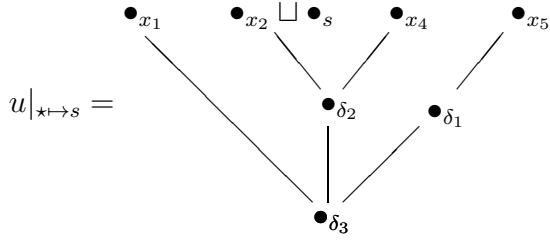
 \setlength{\unitlength}{1cm}
\begin{picture}(4,3.4)
\put(7,0.5){\line(1,1){1}} \put(6.8,0.5){\line(0,1){1}}
\put(6.6,0.5){\line(-1,1){2.2}} \put(6.6,0.2){$\bullet_{\delta_3}$}
\put(6.7,2.0){\line(-4,5){0.6}} \put(7.0,2.0){\line(4,5){0.6}}
\put(8.5,2.0){\line(4,5){0.6}} \put(2.6,1.7){$u|_{\star \mapsto
s}=$}
\put(6.6,0.2){$\bullet_{\delta_3}$}\put(6.7,1.7){$\bullet_{\delta_2}$}
\put(4.1,2.9){$\bullet_{x_1}$}\put(8.1,1.6){$\bullet_{\delta_1}$}
\put(5.5,2.9){$\bullet_{x_2}\sqcup \bullet_{s}$}
\put(7.6,2.9){$\bullet_{x_4}$}\put(9.2,2.9){$\bullet_{x_5}$}
\end{picture}\\

where in $u$ we has ``$\star$" instead of ``$s$".

Similarly, we can define  $(\star_1, \star_2)$-$\Omega$-words as
expressions in  $\mathfrak{S}(X\cup \{\star_1, \star_2\})$  with
only one occurrence of $\star_1$ and only one occurrence of
$\star_2$. Let us denote by $\mathfrak{S}^{\star_1, \star_2} (X)$
the set of all $(\star_1, \star_2)$-$\Omega$-words. Let $u\in
\mathfrak{S}^{\star_1, \star_2} (X)$. Then we call

$$
u|_{s_1,\ s_2}= u|_{\star_1\mapsto s_1,\star_2\mapsto s_2},
$$
an $s_1$-$s_2$-$\Omega$-word.

Now, we assume that $ \mathfrak{S}(X)$ is equipped with a monomial
ordering $>$. This means that $>$ is a well ordering on
$\mathfrak{S}(X)$ such that for any $ v, w \in \mathfrak{S}(X)$,
$u\in \mathfrak{S}^\star (X)$,
$$
w>v\Rightarrow u|_w>u|_v.
$$
 Note that such an order on $\mathfrak{S}(X)$ exists, see order
(\ref{o1}) in the next section.

For every $\Omega$-polynomial $f\in K\langle X;\Omega\rangle $, let
$\bar{f}$  be the leading $\Omega$-word of $f$. If the coefficient
of $\bar{f}$ is $1$, then we call $f$ monic.

\ \

Let $f, g$ be two monic $\Omega$-polynomials. Then, there are two
kinds of $compositions$.
\begin{enumerate}
\item[(I)]If there exists a $\Omega$-word $w=\bar{f}a=b\bar{g}$ for some $a,b\in
\mathfrak{S}(X)$ such that $bre(w)< bre(\bar{f})+bre(\bar{g})$, then
we call $(f,g)_{w}=fa-bg$ the $intersection$ $composition$ of $f$
and $g$ with respect to $w$.
\item[(II)] If there exists a $\Omega$-word $w=\bar{f}=u|_{\bar{g}}$ for some
$u \in \mathfrak{S}^\star(X)$, then we call $(f,g)_{w}=f-u|_{g}$ the
$including$ $composition$ of  $f$ and $g$ with respect to $w$. In
this case transformation $f\mapsto (f,g)_{w}$ is called the
Elimination of the Leading Word (ELW) of $g$ in $f$.
\end{enumerate}

In the above definition, $w$ is called the $ambiguity$ of the
composition. Clearly,
$$
(f,g)_w\in Id(f,g) \ \ \ \mbox{ and }  \ \ \ \overline{(f,g)_w}< w,
$$
where $Id(f,g)$ is the $\Omega$-ideal of $K\langle X;\Omega\rangle$
generated by $f,\ g$.

Let $S$ be a set of monic $\Omega$-polynomials. Then the composition
$(f,g)_w$ is called trivial modulo $(S,w)$ if
$$
(f,g)_w=\sum\alpha_iu_i|_{s_i},
$$
where each $\alpha_i\in K$,  $u_i\in \mathfrak{S}^\star(X)$, $s_i\in
S$ and $u_i|_{\overline{s_i}}< w$. If this is the case, we write
$$
(f,g)_w\equiv 0 \ \ mod (S,w).
$$
In general, for any two $\Omega$-polynomials $p$ and $q$, $ p\equiv
q \ \ mod (S,w) $ means that $ p-q=\sum\alpha_iu_i|_{s_i}, $ where
each $\alpha_i\in K$,  $u_i\in \mathfrak{S}^\star(X)$, $s_i\in S$
and $u_i|_{\overline{s_i}}< w$.

$S$ is called a Gr\"{o}bner-Shirshov basis  in  $K\langle
X;\Omega\rangle$ if any composition $(f,g)_w$ of $f,g\in S$  is
trivial modulo $(S,w)$.

\begin{lemma}\label{l3.3}
Let $S$ be a Gr\"{o}bner-Shirshov  basis  in  $K\langle
X;\Omega\rangle$, $u_1,u_2\in \mathfrak{S}^{\star}(X)$
 and $s_1, s_2\in S$. If
 $w=u_1|_{\overline{s_1}}=u_2|_{\overline{s_2}}$, then
$$
u_1|_{s_1}\equiv u_2|_{s_2} \ \ mod (S,w).
$$
\end{lemma}
{\bf Proof:} \ There are three cases to consider.

(I)\ \ The $\Omega$-words $\overline{s_1}$ and $\overline{s_2}$ are
disjoint in $w$. Then there exits a
$(\star_1,\star_2)$-$\Omega$-word $\Pi$ such that
$$\Pi|_{\overline{s_1},\
\overline{s_2}}=u_1|_{\overline{s_1}}=u_2|_{\overline{s_2}}.
$$
Then
\begin{eqnarray*} u_2|_{ s_2}-u_1|_{
s_1}&=&\Pi|_{\overline{s_1}, \ s_2}-\Pi|_{s_1, \ \overline{s_2}}\\
&=&(-\Pi|_{s_1-\overline{s_1}, \ s_2}+\Pi|_{s_1, \
s_2})+(\Pi|_{s_1,\ s_2-\overline{s_2}}-\Pi|_{s_1, \ s_2})\\
&=&-\Pi|_{s_1-\overline{s_1}, \ s_2}+\Pi|_{s_1,\
s_2-\overline{s_2}}.
\end{eqnarray*}

Let $$-\Pi|_{s_1-\overline{s_1}, \
s_2}=\sum\alpha_{2_t}u_{2_t}|_{s_2}\ \ \mbox{and}\ \ \Pi|_{s_1,\
s_2-\overline{s_2}}=\sum\alpha_{1_l}u_{1_l}|_{s_1}.$$
 Since
$\overline{s_1-\overline{s_1}}<\overline{s_1}$ and
$\overline{s_2-\overline{s_2}}<\overline{s_2}$, we have $
u_{2_t}|_{\overline{s_2}}, \ u_{1_l}|_{\overline{s_1}}<w.$ Therefore
$$
u_2|_{ s_2}-u_1|_{
s_1}=\sum\alpha_{2_t}u_{2_t}|_{s_2}+\sum\alpha_{1_l}u_{1_l}|_{s_1}
$$
with each $u_{2_t}|_{\overline{s_2}}, \
u_{1_l}|_{\overline{s_1}}<w$. It follows that

$$
u_1|_{s_1}\equiv u_2|_{s_2} \ mod (S,w).
$$

(II) The $\Omega$-words  $\overline{s_1}$ and
 $\overline{s_2}$ have nonempty intersection in $w$ but do not
 include each other. Without lost of generality we can assume that
$ \overline{s_1}a=b\overline{s_2}$ is a $\Omega$-subword of $w$ for
some $\Omega$-words $a,\ b$. Then there exists a
$\star$-$\Omega$-word $\Pi$ such that
$$
\Pi|_{\overline{s_1}a}=u_1|_{\overline{s_{_1}}}=u_2|_{\overline{s_2}}=\Pi|_{b\overline{s_2}}.
$$
Thus, we have
$$
u_{_2}|_{ s_{_2}}-u_1|_{ s_{_1}}=\Pi|_{bs_{_2}}-\Pi|_{s_{_1}a}
=-\Pi|_{s_{_1}a-bs_{_2}}.
$$
Since  $S $ is a Gr\"{o}bner-Shirshov basis in  $K\langle
X;\Omega\rangle$,
 we have
$$
s_1a-bs_2=\sum\alpha_jv_j|_{s_j},
$$
where each $\alpha_j\in K, \ v_j\in \mathfrak{S}^\star(X), \ s_j\in
S$ and $v_j|_{\overline{s_j}}< \overline{s_1}a$. Let $
\Pi|_{v_j|_{s_j}}=\Pi_{j}|_{s_j}. $ Then $ u_2|_{ s_2}-u_1|_{
s_1}=\sum\alpha_j\Pi_{j}|_{s_j} $ with $
\Pi_{j}|_{\overline{s_j}}<w. $ It follows that
$$
u_1|_{s_1}\equiv u_2|_{s_2} \ \  mod (S,w).
$$

(III) One of $\Omega$-words $\overline{s_1}$,
 $\overline{s_2}$ is contained in the other.
For example, let $ \overline{s_1}=u|_{\overline{s_2}}$ for some
$\star$-$\Omega$-word $u$. Then
$$
w=u_2|_{\overline{s_2}}=u_1|_{u|_{\overline{s_2}}}\ \ \mbox{ and }\
\ u_2|_{ s_2}-u_1|_{ s_1}=u_1|_{u|_{s_2}}-u_1|_{ s_1}=-u_1|_{
s_1-u|_{s_2}}.
$$
By using similar method as in (II), we  obtain the result. \hfill
$\blacksquare$
\\

The following theorem is an analogy  of  Shirshov's Composition
lemma for Lie algebras \cite{Sh62b}, which was specialized to
associative algebras by Bokut \cite{Bo76}, see also Bergman
\cite{Be78}.

\begin{theorem}\label{3.5}{\em(Composition-Diamond lemma)}\ \  Let $S$ be a set of monic
$\Omega$-polynomials in $K\langle X;\Omega\rangle$,
 $>$ a monomial ordering on $\mathfrak{S}(X)$ and $Id(S)$ the $\Omega$-ideal of
 $K\langle X;\Omega\rangle$ generated by $S$.  Then the following
statements are equivalent:
 \begin{enumerate}
\item[(I)] $S $ is a Gr\"{o}bner-Shirshov basis in $K\langle X;\Omega\rangle$.
\item[(II)] $ f\in Id(S)\Rightarrow \bar{f}=u|_{\overline{s}}$
for some $u \in \mathfrak{S}^\star(X)$ and $s\in S$.
\item[$(II')$]  $f\in Id(S)\Rightarrow
f=\alpha_1u_1|_{s_1}+\alpha_2u_2|_{s_2}+\cdots+\alpha_nu_n|_{s_n}$
where each $\alpha_i\in K, \ s_i\in S,\ u_i\in
\mathfrak{S}^\star(X)$ and
$u_1|_{\overline{s_1}}>u_2|_{\overline{s_2}}>\cdots>u_n|_{\overline{s_n}}$.
\item[(III)] $Irr(S) = \{ w\in \mathfrak{S}(X) |  w \neq
u|_{\overline{s}}
 \mbox{ for  any} \ u \in \mathfrak{S}^\star(X) \ \mbox{and } s\in S\}$
is a $K$-basis of $K\langle X;\Omega|S\rangle=K\langle
X;\Omega\rangle/Id(S)$.
\end{enumerate}
\end{theorem}
{\bf Proof:}\ \

 (I)$\Longrightarrow$ (II)\ \ Let  $0\neq f\in
Id(S)$. Then

$$
f=\sum\limits_{i=1}^{n}\alpha_i u_i|_{ s_i}
$$
where each $\alpha_i\in K$, $u_i\in \mathfrak{S}^\star(X)$ and $
s_i\in S $.

Let $w_i=u_i|_{\overline{ s_i}}$. We arrange these leading
$\Omega$-words in non-increasing order by
$$
w_1= w_2=\cdots=w_m >w_{m+1}\geq \cdots\geq w_n.
$$
Now we prove the result by  induction  on $m$.

If $m=1$, then $\bar{f}=u_1|_{\overline{ s_1}}$.

Now we assume that $m\geq 2$. Then $
u_1|_{\overline{s_1}}=w_1=w_2=u_2|_{\overline{s_2}}$.

We prove the result by induction on $w_1$.  It is clear that
$w_1\geq \overline{f}$. If $w_1= \overline{f}$ then the result is
clear. Let $w_1> \overline{f}$. Since $S $ is a Gr\"{o}bner-Shirshov
basis in $K\langle X;\Omega\rangle$, by Lemma \ref{l3.3}, we have
$$
u_2|_{ s_2}-u_1|_{ s_1}=\sum\beta_jv_j|_{s_j}
$$
where $\beta_j\in K, \ s_j\in S, v_j\in \mathfrak{S}^\star(X)$  and
$v_j|_{\overline{s_j}}<w_1$. Therefore, since
$$
 \alpha_1u_1|_{ s_1}+\alpha_2u_2|_{
s_2}=(\alpha_1+\alpha_2)u_1|_{ s_1}+\alpha_2(u_2|_{ s_2}-u_1|_{
s_1}),
$$
we have
$$
f=(\alpha_1+\alpha_2)u_1|_{ s_1}+\sum\alpha_2\beta_jv_j|_{s_j}+
\sum\limits_{i=3}^{n}\alpha_iu_i|_{ s_i}.
$$

If either $m>2$ or $\alpha_1+\alpha_2\neq 0$, then the result
follows from induction on $m$. If $m=2$ and $\alpha_1+\alpha_2=0$,
then  the result follows from induction on $w_1$.

From the above proof, it follows also $(II)\Longleftrightarrow
(II')$.

(II)$\Longrightarrow$ (III) For any $f\in K\langle X;\Omega\rangle$,
by induction on $\bar{f}$, using the ELW's of $S$, we have
$$
f=\sum \limits_{u_i\in Irr(S),u_i\leq\bar{h}}\alpha_iu_i + \sum
\limits_{s_j\in S,
v_j|_{\overline{s_j}}\leq\bar{h}}\beta_jv_j|_{s_j}.
$$

It follows that any $f\in K\langle X;\Omega\rangle$
 can be expressed modulo $Id(S)$ by a linear combination of elements of $Irr(S)$.
 Now suppose that $\alpha_1u_1+\alpha_2u_2+\cdots\alpha_nu_n=0$ in
$ K\langle X;\Omega|S\rangle$ with each $u_i\in Irr(S)$.
  Then, in $ K\langle
X;\Omega\rangle$,
$$
g=\alpha_1u_1+\alpha_2u_2+\cdots+\alpha_nu_n\in Id(S).
$$
By (II), we have $\bar{g}\notin Irr(S)$, a contradiction.  This
shows that $Irr(S)$ is a $K$-basis of $ K\langle X;\Omega|S\rangle$.

(III)$\Longrightarrow $(I) In particular, for any composition
$(f,g)_w$ in $S$, since $(f,g)_w\in Id(S)$ and (III), we have
$$
(f,g)_w=\sum\beta_jv_j|_{s_j},
$$
where each $\beta_j\in K$,  $v_j\in \mathfrak{S}^\star(X)$, $s_j\in
S$ and $v_j|_{\overline{s_j}}\leq\overline{(f,g)_w}< w$.  This shows
(I). \hfill $\blacksquare$

\ \

\noindent{\bf Remark:}   If  $\Omega$ is empty, then  Theorem
\ref{3.5} is the  Composition-Diamond Lemma for associative algebras
(cf. \cite{Bo76}).

\section{ Gr\"{o}bner-Shirshov bases for  free Rota-Baxter algebras}

In this section, we give a Gr\"{o}bner-Shirshov basis and a linear
basis for a free Rota-Baxter algebra. In fact, we construct the free
Rota-Baxter algebra on a set $X$ by using the Composition-Diamond
Lemma  above for associative algebras with multiple linear operators
(Theorem \ref{3.5}).

First of all,  we define an order on $\mathfrak{S}(X)$, which will
be used in the sequel. Let $X$ and $\Omega$ be well ordered sets. We
order $X^*$ by the deg-lex ordering. For any $u\in \mathfrak{S}(X)$,
$u$ can be uniquely expressed without brackets as
$$
u=u_0\delta_{i_{_{1}}}\overrightarrow{x_{i_1}}u_1\cdots
\delta_{i_{_{t}}}\overrightarrow{x_{i_{t}}}u_t,
$$
where each $u_i\in X^*,\delta_{i_{_{k}}}\in \Omega_{i_{_{k}}}, \
\overrightarrow{x_{i_k}}=( x_{k_{1}},x_{k_{2}},\cdots,
x_{k_{i_{_{k}}}} )\in \mathfrak{S}(X)^{i_k}$. It is reminded that
for each $i_{{k}}, \ dep(u)>dep(\overrightarrow{x_{i_k}})$.

Denote by
$$
wt(u)=(t,\delta_{i_{_{1}}},\overrightarrow{x_{i_{1}}},
\cdots,\delta_{i_{_{t}}},\overrightarrow{x_{i_t}}, u_0, u_1, \cdots,
u_t ).
$$
Let $v=v_0\mu_{j_{_{1}}}\overrightarrow{y_{j_{1}}}v_1\cdots
\mu_{j_{_{s}}}\overrightarrow{y_{j_{s}}}v_s\in  \mathfrak{S}(X). $
Then, we order  $\mathfrak{S}(X)$ as follows:
\begin{equation}\label{o1}
u>v\Longleftrightarrow wt(u)>wt(v)\ \mbox{ lexicographically}
\end{equation}
by induction on $dep(u)+dep(v)$.

It is clear that the order (\ref{o1}) is a monomial ordering on
$\mathfrak{S}(X)$.

\ \

 Let $K$ be a commutative ring with unit and
$\lambda\in K$. A Rota-Baxter $K$-algebra of weight $\lambda$
(\cite{Bax60, Ro69}) is an associative algebra $R$ with  a
$K$-linear operation $ P:R\rightarrow R$ satisfying the Rota-Baxter
relation:
$$
 P(x)P(y) = P( P(x)y +xP(y)+\lambda xy), \  \forall x,y \in  R.
$$

Thus, any Rota-Baxter algebra is a special case of associative
algebra with multiple operators when  $\Omega=\{P\}$.

Now, let  $\Omega=\{P\}$ and $\mathfrak{S}(X)$ be as before. Let
$K\langle X;P\rangle$ be the free associative algebra with one
operator $\Omega=\{P\}$ on a set $X$.

\begin{theorem}\label{t4.1}
With the order (\ref{o1})  on $\mathfrak{S}(X)$,
$$
S=\{P(x)P(y) - P( P(x)y) - P(xP(y))-\lambda P(xy) |\  x,y \in
\mathfrak{S}(X)\}
$$
is a Gr\"{o}bner-Shirshov basis in $K\langle X;P\rangle$.
\end{theorem}
{\bf Proof:} The ambiguities of  all possible compositions of
$\Omega$-polynomials in $S$ are:
$$
(i) \ \ \ \  P(x)P(y)P(z) \ \ \ \ (ii) \ \ \ \ P(u|_{P(x)P(y)})P(z)\
\ \ \ (iii) \ \ \ \  P(x)P(u|_{P(y)P(z)})
$$
where $x, y,z\in \mathfrak{S}(X)$, $u \in \mathfrak{S}^\star(X)$. It
is easy to check that all these compositions are trivial. Here, for
example, we just check $(ii)$. Others are similarly proved. Let
$$
f(x,y)=P(x)P(y) - P( P(x)y) - P(xP(y))-\lambda P(xy).
$$
Then
\begin{eqnarray*}
&&(f(u|_{P(x)P(y)}, z), f(x,y))_{_{P(u|_{P(x)P(y)})P(z)}}\\
&&=-P(P(u|_{P(x)P(y)})z)-P(u|_{P(x)P(y)}P(z))-\lambda P(u|_{P(x)P(y)}z)\\
&& \ \  \  +P(u|_{P( P(x)y)})P(z)+P(u|_{P( xP(y))})P(z)+\lambda P(u|_{P(xy)})P(z)\\
&&\equiv -P(P(u|_{P( P(x)y)})z)-P(P(u|_{P( xP(y))})z)-\lambda P(P(u|_{ P(xy)})z)\\
&& \ \  \ -P(u|_{P( P(x)y)}P(z))-P(u|_{P( xP(y))}P(z))-\lambda P(u|_{P( xy)}P(z))\\
&& \ \  \ -\lambda P(u|_{P(P(x)y)}z)-\lambda P(u|_{P(xP(y))}z)-\lambda^2 P(u|_{p(xy)}z)\\
&& \ \  \ +P(P(u|_{P( P(x)y)})z) + P(u|_{P( P(x)y)}P(z))+\lambda P(u|_{P( P(x)y)}z)\\
&& \ \  \ +P(P(u|_{P( xP(y))})z)+P(u|_{P( xP(y))}P(z))+\lambda P(u|_{P( xP(y))}z)\\
&& \ \  \ +\lambda P(P(u|_{P(xy)})z)+\lambda P(u|_{P(xy)}P(z))+\lambda^2P(u|_{P(xy)}z)\\
&&\equiv 0 \ \ mod(S,P(u|_{P(x)P(y)})P(z)). \ \ \ \ \ \ \ \ \ \ \ \
\ \ \ \  \ \ \ \ \ \ \ \  \  \  \ \   \ \ \ \ \ \ \ \ \ \ \ \ \
\blacksquare
\end{eqnarray*}

Now, we describe the set $Irr(S)$. Let $Y$ and $Z$ be two subsets of
$\mathfrak{S}(X)$. Define the alternating product of $X$ and $Y$
with  $P$ (see also \cite{EG08a}):
$$
\Lambda_{X}^{P}(Y,Z)=(\cup_{ r\geq 1}(YP(Z))^r)\cup (\cup_{r\geq
0}(YP(Z))^rY) \cup(\cup_{ r\geq
1}(P(Z)Y)^r)\cup(\cup_{r\geq0}(P(Z)Y)^rP(Z)).
$$
Clearly, $\Lambda_{X}^{P}(Y,Z)\subseteq \mathfrak{S}(X)$. Now let
$$
\Phi_0=S(X).
$$
For $n>0$, define
$$
\Phi_n=\Lambda_{X}^{P}(\Phi_0, \Phi_{n-1}).
$$
Let
$$
\Phi(X)=\bigcup_{n\geq 0 }\Phi_n.
$$
Then, we call the elements in $\Phi(X)$ the Rota-Baxter  words. It
is easy to see that
$$
\Phi_0\subset\Phi_1\subset\cdots\subset\Phi_n\subset\cdots
$$

\begin{theorem}\label{t4.3}{\em (\cite {EG08a})}
$Irr(S)= \Phi(X)$  is a linear basis of  $K\langle X;P|S\rangle$.
\end{theorem}
{\bf Proof:} The result follows from Theorem \ref{3.5} and Theorem
\ref{t4.1}.  \hfill $\blacksquare$

Thus, we have the following theorem.

\begin{theorem} \label{t4.3}{\em (\cite {EG08a})}
$K\langle X;P|S\rangle$ is a free Rota-Baxter algebra on  set $X$
with a $K$-basis $\Phi(X)$.
\end{theorem}

  By using ELWs, we have the following algorithm. In fact, it is    an  algorithm to
  compute the product of two  Rota-Baxter words in the free Rota-Baxter algebra $K\langle X;P|S\rangle$.
It coincides with the product defined in \cite{EG08a}.

\ \

 \noindent
\textbf{Algorithm} Let $u,v\in  \Phi(X)$. We define $u\diamond v$ by
induction on $n=dep(u)+dep(v)$.

 (a) If $n=0$, then $u,v \in S(X)$ and $u\diamond
v=uv$.

 (b)  If $n\geq 1$, there are two cases to consider:

 (i) If  $bre(u)=bre(v)=1$, i.e., $u, v\in X\cup
 P(\Phi(X))$, then
 \newline
\begin{equation*}
 u\diamond v=\left\{
\begin{array}{l@{\quad\quad}l}
uv & \ \mbox{if}\  u\in X \  \mbox{or} \ v\in X  \\
P(P(u^{'})\diamond v^{'}) + P(u^{'}\diamond P(v^{'}))+\lambda
P(u^{'}\diamond v^{'}) & \
\mbox{if} \  u=P(u^{'})\  \mbox{and} \ v=P(v^{'})%
\end{array}%
\right.
\end{equation*}

(ii) If $bre(u)>1$ or $bre(v)>1$, say, $ u=u_1u_2\cdots u_t \ \ \
\mbox{ or }\ \ \ v=v_1v_2\cdots v_l, $ where $u_i$ and $v_j$ are
prime, then
$$
 u\diamond v =u_1u_2\cdots u_{t-1} (u_t\diamond v_1)v_2\cdots v_l.
$$

\section{Gr\"{o}bner-Shirshov bases for free $\lambda$-differential
algebras}

In this section, we give two Gr\"{o}bner-Shirshov bases and two
linear bases for a free $\lambda$-differential algebra. We construct
 the free $\lambda$-differential algebra on a set $X$ by
using the Composition-Diamond Lemma (Theorem \ref{3.5}).

Let $K$ be a commutative ring with unit and $\lambda\in K$. A
$\lambda$-differential
 algebra over $K$ (\cite{GK08})  is an  associative $K$-algebra  $R$
  together with a  $K$-linear operator $D:R\rightarrow R$ such that
$$
D(xy)=D(x)y+xD(y)+\lambda D(x)D(y),\ \forall x, y \in R.
$$

Any  $\lambda$-differential  algebra is also an   associative
algebra with one
 operator $\Omega=\{D\}$.

Let $X$ be well ordered and $K\langle X;D\rangle$ the free
associative algebra with one operator  $\Omega=\{D\}$ defined as
before.

 For any $u\in
\mathfrak{S}(X)$, $u$ has a unique expression
$$
u=u_1u_2\cdots u_n,
$$
where each $u_i\in X\cup D(\mathfrak{S}(X))$. Denote by
$deg_{_{X}}(u)$ the number of $x\in X$ in $u$. For example, if
$u=D(x_1x_2)D(D(x_1))x_3\in \mathfrak{S}(X)$, then
$deg_{_{X}}(u)=4$.  Let
$$
wt(u)=(deg_{_{X}}(u), u_1,u_2,\cdots, u_n).
$$
Let  $v=v_1v_2\cdots v_m\in \mathfrak{S}(X)$. Define
\begin{equation}\label{o2}
u>v\Longleftrightarrow wt(u)>wt(v)\ \mbox{ lexicographically},
\end{equation}
where for each  $t, \ u_t>v_t$ if one of the following holds:

(a) $u_t, v_t\in X$ and $u_t>v_t$;

(b)  $u_t=D(u_{t}^{'}), v_t\in X$;

(c)  $u_t=D(u_{t}^{'}),v_t=D(v_{t}^{'})$ and $u_{t}^{'}>v_{t}^{'}$.

\ \

It is easy to see the order (\ref{o2}) is a monomial ordering on
$\mathfrak{S}(X)$.

\begin{theorem}\label{t5.3}With  the  order
(\ref{o2})  on $\mathfrak{S}(X)$,
$$
S=\{D(xy)-D(x)y-xD(y)-\lambda D(x)D(y) |\  x,y \in \mathfrak{S}(X)\}
$$
is a Gr\"{o}bner-Shirshov basis in $K\langle X;D\rangle$.
\end{theorem}
{\bf Proof:} The ambiguities of  all possible compositions of the
$\Omega$-polynomials in $S$ are:
$$
(a)\ \ \ \ \ D(v(u|_{_{D(xy)}})) \ \ \ \ \ \ \ (b) \ \ \ \
D(u|_{_{D(xy)}}v)
$$
where $x, y, v\in \mathfrak{S}(X), u\in \mathfrak{S}^\star(X)$. It
is easy to check that all these compositions are trivial. Here, for
example, we just check $(b)$. Others are similarly proved. Let
$$
g(x,y)=D(xy)-D(x)y-xD(y)-\lambda D(x)D(y).
$$

Then
\begin{eqnarray*}
&&(g(u|_{_{D(xy)}}, v),g(x,y))_{_{D((u|_{_{D(xy)}})v)}}\\
&&=-D(u|_{_{D(xy)}})v-u|_{_{D(xy)}}D(v)-\lambda
D(u|_{_{D(xy)}})D(v)\\
&& \ \  \  +D((u|_{_{D(x)y}})v)+D((u|_{_{xD(y)}})v)+\lambda D((u|_{_{D(x)D(y)}})v) \\
&&\equiv -D(u|_{_{D(x)y}})v-D(u|_{_{xD(y)}})v-\lambda D(u|_{_{D(x)D(y)})}v\\
&& \ \  \ -u|_{_{D(x)y}}D(v)-u|_{_{xD(y)}}D(v)-\lambda u|_{_{D(x)D(y)}}D(v)\\
&& \ \  \ -\lambda D(u|_{_{D(x)y}})D(v)-\lambda D(u|_{_{xD(y)}})D(v)-\lambda^2 D(u|_{_{D(x)D(y)}})D(v)\\
&& \ \  \ +D(u|_{_{D(x)y}})v + u|_{_{D(x)y}}D(v)+\lambda D(u|_{_{D(x)y}})D(v)\\
&& \ \  \ +D(u|_{_{xD(y)}})v+u|_{_{xD(y)}}D(v)+\lambda D(u|_{_{xD(y)}})D(v)\\
&& \ \  \ +\lambda D(u|_{_{D(x)D(y)}})v+\lambda  u|_{_{D(x)D(y)}}D(v)+\lambda^2 D(u|_{_{D(x)D(y)}})D(v)\\
&&\equiv 0 \ \ mod(S,D((u|_{_{D(xy)}})v)). \ \ \ \ \ \ \ \ \ \ \ \ \
\ \ \  \ \ \ \ \ \ \ \  \  \  \ \   \ \ \ \ \ \ \ \ \ \ \ \ \
\blacksquare
\end{eqnarray*}

 Let $D^{\omega}(X)=\{D^i(x)|i\geq
0, x\in X\}$, where $D^0(x)=x$  and $S(D^{\omega}(X))$ the free
semigroup generated by $D^{\omega}(X)$.

\ \

\begin{theorem}\label{t5.4}{\em(\cite{GK08})}
$Irr(S)= S(D^{\omega}(X))$  is a $K$-basis of $K\langle
X;D|S\rangle$.
\end{theorem}
{\bf Proof:} The result follows from Theorem \ref{3.5} and Theorem
\ref{t5.3}.  \hfill $\blacksquare$

Thus, we have the following theorem.

\begin{theorem}\label{t5.5}{\em(\cite{GK08})}
$K\langle X;D|S\rangle$ is a free $\lambda$-differential  algebra on
set $X$ with a linear basis $S(D^{\omega}(X))$.
\end{theorem}

By using ELWs, we have the following algorithm.

\noindent  \textbf{Algorithm} (\cite{GK08}) Let $u=u_1u_2\cdots u_t
\in S(D^{\omega}(X))$, where $u_i\in D^{\omega}(X)$.     Define
$D(u)$ by induction on $t$.

 (a) If $t=1$, i.e., $u=D^i(x)$ for some $i\geq 0, x\in X $,
 then $D(u)=D^{(i+1)}(x)$.

 (b) If $t\geq 1$, then
$$
D(u)=D(u_1u_2\cdots u_t )=D(u_1)(u_2\cdots u_t )+u_1D(u_2\cdots u_t)
+\lambda D(u_1)D(u_2\cdots u_t).
$$

\ \

The following theorem gives another Gr\"{o}bner-Shirshov basis in
$K\langle X;D\rangle$ when  $\lambda\in K$ is invertible.

\begin{theorem}\label{t5.6} Let  $\lambda\in K$ be invertible
and $ \mathfrak{S}(X)$ with the  order (\ref{o1}). Set
$$
 T=\{D(x)D(y)+\lambda^{-1}D(x)y+\lambda^{-1}xD(y)-\lambda^{-1}D(xy)|x, y\in
\mathfrak{S}(X)\}.
$$
Then $K\langle X;D|S\rangle=K\langle X;D|T\rangle$ and   $T$ is a
Gr\"{o}bner-Shirshov basis in $K\langle X;D\rangle$.
\end{theorem}
{\bf Proof:} The proof is similar to the Theorem \ref{t4.1}. \hfill
$\blacksquare$

\ \

 Now let
$$
\Phi_0=S(X).
$$
For $n>0$, define
$$
\Phi_n=\Lambda_{X}^{D}(\Phi_0, \Phi_{n-1}).
$$
Let
$$
\Phi(X)=\bigcup_{n\geq 0 }\Phi_n.
$$

\begin{theorem}\label{t5.7} Let  $\lambda\in K$ be invertible. Then
$Irr(T)= \Phi(X)$  is a $K$-basis of $K\langle X;D|T\rangle$.
\end{theorem}
{\bf Proof:} By  Theorem \ref{t5.6} and Theorem \ref{3.5}, we can
obtain the result. \hfill $\blacksquare$

By  Theorem \ref{t5.7} and  Theorem \ref{t5.5}, we have the
following theorem.
\begin{theorem} Let  $\lambda\in K$ be invertible. Then
$K\langle X;D|T\rangle$ is a free $\lambda$-differential algebra on
set $X$ with a  linear basis $\Phi(X)$ as above.
\end{theorem}

\section{ Gr\"{o}bner-Shirshov bases for free  $\lambda$-differential Rota-Baxter  algebras}

In this section, we obtain a Gr\"{o}bner-Shirshov basis and a linear
basis for a free $\lambda$-differential Rota-Baxter algebra.

Let $K$ be a commutative ring with unit and $\lambda\in K$.
 A differential Rota-Baxter algebra of weight $\lambda$
 (\cite{GK08}), called also $\lambda$-differential Rota-Baxter
 algebra,
is an associative $K$-algebra $R$ with two $K$-linear operators
$P,D:R\rightarrow R$ such that for any $x,y\in R$,

(I) (Rota-Baxter relation) $P(x)P(y)=P(xP(y))+P(P(x)y)+\lambda
P(xy);$

(II) ($\lambda$-differential relation) $ D(xy)=D(x)y+xD(y)+\lambda
D(x)D(y);$

(III) $D(P(x))=x$.

Hence, any  $\lambda$-differential Rota-Baxter
 algebra is an    associative algebra with  two linear operators $\Omega=\{P,
 D\}$.

Let  $K\langle X;\Omega \rangle$ be  the free  associative algebra
with multiple linear operators $\Omega$ on $X$, where
$\Omega=\{P,D\}$. For any $u\in \mathfrak{S}(X)$, $u$ has a unique
expression
$$
u=u_0P(b_1)u_1P(b_2)u_2\cdots P(b_n) u_n,
$$
where each $u_i\in (X\cup D(\mathfrak{S}(X)))^*$ and $b_i\in
\mathfrak{S}(X)$. Denote by
$$
wt(u)=(deg_{_{X}}(u),deg_{_{P}}(u),n, b_1, \cdots, b_n,
u_0,\cdots,u_n),
$$
where $deg_{_{P}}(u)$ is the number of $P$ in $u$. Also, for any
$u_t\in (X\cup D(\mathfrak{S}(X)))^* $, $u_t$ has a unique
expression
$$
u_t=u_{t_{1}}\cdots  u_{t_{k}}
$$
where each $u_{t_{j}}\in X\cup D(\mathfrak{S}(X))$.

 Let $X$  be well ordered  and $v=v_0P(c_1)v_1P(c_2)v_2\cdots
P(c_m) v_m$. Order  $\mathfrak{S}(X)$ as follows:
\begin{equation}\label{o3}
u>v\Longleftrightarrow wt(u)>wt(v)\ \mbox{ lexicographically}
\end{equation}
 where for each $t$, $\ u_t>v_t$ if
 $$
 (deg_{_{X}}(u_t),deg_{_{P}}(u_t), u_{t_{1}},\cdots,u_{t_{k}})>
 (deg_{_{X}}(v_t),deg_{_{P}}(v_t), v_{t_{1}},\cdots,v_{t_{l}})\ \mbox{
 lexicographically}
 $$
where for each $j$,  $u_{t_{j}}>v_{t_{j}}$ if  one of the following
holds:

(a)  $u_{t_{j}},v_{t_{j}}\in X$ and $u_{t_{j}}>v_{t_{j}}$;

(b)  $u_{t_{j}}=D(u_{t_{j}}^{'}), v_{t_{j}}\in X$;

(c)  $u_{t_{j}}=D(u_{t_{j}}^{'}),v_{t_{j}}=D(v_{t_{j}}^{'})$ and
$u_{t_{j}}^{'}>v_{t_{j}}^{'}$.

\ \

Clearly, the order (\ref{o3}) is a monomial ordering on
$\mathfrak{S}(X)$.

Let $S$ be the set consisting of  the following
$\Omega$-polynomials:
\begin{enumerate}
\item[1.]
$P(x)P(y)-P(xP(y))-P(P(x)y)-\lambda P(xy)$,
\item[2.] $D(xy)-D(x)y-xD(y)-\lambda D(x)D(y)$,
\item[3.] $D(P(x))-x$,
\end{enumerate}
where $x,y\in \mathfrak{S}(X)$.

\begin{theorem}\label{t6.3}With  the order
(\ref{o3}) on $\mathfrak{S}(X)$,  $S$ is a Gr\"{o}bner-Shirshov
basis in $K\langle X;\Omega\rangle$.
\end{theorem}
\noindent {\bf Proof.}  Denote by $i\wedge j$ the composition of the
$\Omega$-polynomials of type $i$ and type $j$. The ambiguities of
all possible compositions of the $\Omega$-polynomials  in $S$
 are only as below:
 \begin{tabbing}
 $3\wedge3\ \ \ D(P(u|_{_{D(P(x))}}))$\\[0.7ex]
 $3\wedge2\ \ \ D(P(u|_{_{D(xy)}}))$\\[0.7ex]
 $3\wedge1\ \ \ D(P(u|_{_{P(x)P(y)}}))$\\[0.7ex]
  $2\wedge3\ \ \ D(x(u|_{_{D(P(y))}})),\ \ \
  D((u|_{_{D(P(x))}})y)$\\[0.7ex]
  $2\wedge2\ \ \ D(x(u|_{_{D(yz)}})),\ \ \
  D((u|_{_{D(xy)}})z)$\\[0.7ex]
  $2\wedge1\ \ \ D((u|_{_{P(x)P(y)}})z),\ \ \
  D(x(u|_{_{P(y)P(z)}}))$\\[0.7ex]
   $1\wedge3\ \ \ P(u|_{_{D(P(x))}})P(z),\ \ \
 P(x)P(u|_{_{D(P(y))}})$\\[0.7ex]
 $1\wedge2\ \ \ P(u|_{_{D(xy)}})P(z),\ \ \
 P(x)P(u|_{_{D(yz)}})$\\[0.7ex]
 $1\wedge1\ \ \ P(x)P(y)P(z), \ \  P(u|_{_{P(x)P(y)}})P(z), \ \ \
 P(x)P(u|_{_{P(y)P(z)}})$
\end{tabbing}
where $x,y,z \in \mathfrak{S}(X)$, $u\in \mathfrak{S}^\star(X)$.
Similar to the proofs of Theorem \ref{t4.1} and Theorem \ref{t5.3},
we have $1\wedge1$ and $2\wedge2$ are trivial. Others are also
easily checked. Here we just check $3\wedge2$:
\begin{eqnarray*}
3\wedge2&=& -u|_{_{D(xy)}}+ D(P(u|_{_{D(x)y}}))+D(P(u|_{_{xD(y)}}))+\lambda D(P(u|_{_{D(x)D(y)}}))\\
&\equiv& -u|_{_{D(x)y}}-u|_{_{xD(y)}}-\lambda u|_{_{D(x)D(y)}} +u|_{_{D(x)y}}+u_{_{xD(y)}}+\lambda u|_{_{D(x)D(y)}}\\
&\equiv& 0 \ \ mod (S,D(P(u|_{_{D(xy)}}))).
\end{eqnarray*}\hfill $\blacksquare$

 Let $D^{\omega}(X)=\{D^i(x)|i\geq
0, x\in X\}$, where $D^0(x)=x$. Define
$$
\Phi_0=S(D^{\omega}(X)),
$$
and for $n>0$,
$$
\Phi_n=\Lambda_{_{D^{\omega}(X)}}^{^{P}}(\Phi_0, \Phi_{n-1}).
$$
Let
$$
\Phi(D^{\omega}(X))=\bigcup_{n\geq 0 }\Phi_n.
$$

\begin{theorem}\label{t6.4}
$Irr(S)= \Phi(D^{\omega}(X))$  is a $K$-basis of  $K\langle
X;\Omega|S\rangle$.
\end{theorem}
{\bf Proof:} By   Theorem \ref{t6.4} and Theorem \ref{3.5}, we can
obtain the result. \hfill $\blacksquare$

Now, we have the following theorem.

\begin{theorem}
$K\langle X;\Omega|S\rangle$ is a free $\lambda$-differential
Rota-Baxter algebra on a set $X$ with a linear basis
$\Phi(D^{\omega}(X))$.
\end{theorem}

In fact, $K\langle X;\Omega|S\rangle$ is a free Rota-Baxter algebra
on $D^{\omega}(X)$ and  this result is similar to the one proved by
L. Guo and W. Keigher in \cite{GK08} for the free
$\lambda$-differential Rota-Baxter algebra with unit  by using the
planar decorated rooted trees.

\ \

\end{document}